\documentclass[10pt,a4paper,twoside]{amsart}
\usepackage[latin1]{inputenc}
\usepackage[T1]{fontenc}

\theoremstyle{plain}
\newtheorem{thm}{\bfseries Theorem}[section]

\newtheorem{prop}[thm]{\bfseries Proposition}

\newtheorem{cor}[thm]{\bfseries Corollary}
\newtheorem{df}[thm]{\bfseries Definition}
\newtheorem{claim}[thm]{\bfseries Claim}
\theoremstyle{remark}

\newtheorem{rem}[thm]{\bfseries Remark}

\numberwithin{equation}{section}
\DeclareMathSymbol{\Z}{\mathalpha}{AMSb}{"5A} 
\DeclareMathSymbol{\PP}{\mathalpha}{AMSb}{"50} 
\DeclareMathSymbol{\Q}{\mathalpha}{AMSb}{"51}
\DeclareMathSymbol{\N}{\mathalpha}{AMSb}{"4E}
\DeclareMathSymbol{\R}{\mathalpha}{AMSb}{"52}

\newcommand{\C} {\mathbb{C}}

\newcommand{\CL} {\mathcal{L}}

\newcommand{\eps}{\varepsilon}

\newcommand{\Wedge}{\bigwedge}

\def \newline{{\par\noindent}}
\def \BZ{{\Z}}
\def \BP{{\PP}}
\def \BN{{\N}}
\def \Sig{{\mathcal S}}
\def \L{{\CL}}
\def \BC{{\C}}
\def \Li{{\rm Li}}
\def \Sym{{\rm Sym}}

\def \BQ{{\Q}}
\def \li{{\ \{}}
\def \re{{\}\ }}
\def \al{{\alpha}}
\def \be{{\beta}}
\def\Sym{{\rm Sym}}

\def \crr {{\rm cr}}
\def \ms {{\medskip}}

\begin{document}
\date{}
\title{Functional equations for higher logarithms}
\author{Herbert Gangl}
\address{MPI Bonn, Vivatsgasse 7, D-53111 Bonn, Germany}
\email{herbert@mpim-bonn.mpg.de}
\maketitle
\setcounter{tocdepth}{1}
\begin{abstract}
Following earlier work by Abel and others, 
Kummer gave functional equations for the polylogarithm
function $\Li_m(z)$ up to $m=5$ in 1850, but no example for larger $m$ was known until recently.
We give the first genuine 2-variable functional equation for the 7--logarithm.
We investigate and relate identities for the 3-logarithm given by 
Goncharov and Wojtkowiak  and deduce a certain family of functional equations for the 4-logarithm.
\end{abstract}

\section{Introduction}
\subsection{}
An essential property of the logarithm is its functional equation 
$$\log(xy)=\log(x)+\log(y)\,.$$
An essential property of the dilogarithm, which is defined by the
power series $$\Li_2(z)=\sum_{n=1}^\infty{}\frac{z^n}{ n^2}\quad, \qquad |z|<1,$$
and can be analytically continued to $\,\BC-[1,\infty)\,$,
is the 5-term relation \big(note the involutary symmetry $(x,y)\mapsto(\frac{1-x}{ 1-y^{-1}},
\frac{1-y}{1-x^{-1}})$\big)
$$\Li_2(xy) -\Li_2(x)-\Li_2(y) -\Li_2\Big(\frac{1-x}{ 1-y^{-1}}\Big)-\Li_2\Big(\frac{1-y}{
1-x^{-1}}\Big)=({\rm  elementary})\,,$$
which was found---in a lot of different forms---by Spence, Abel and many others (cf. \cite{L1}, Chapter I). 
Here the ``elementary'' right hand side consists a sum of products of logarithms. 

In 1840, Kummer \cite{Kummer} found functional equations for higher (poly-)logarithms
$$\Li_m(z)=\sum_{n=1}^\infty{}\frac{z^n}{ n^m}\,,\qquad |z|<1,$$
(which can be likewise analytically continued)
up to $m\le 5$. His method does not extend to higher $m$, though---for a more detailed 
statement we refer the reader to Wechsung's paper \cite{We1}, where more functional equations 
in Kummer's spirit are derived and where it is shown that there are no such ``Kummer-type''
equations for $m>5$.

Other efforts in the direction of finding new equations, even with the help of a computer,
seemed also to be restricted to $m\le 5$ (cf., e.g., \cite{L3}, \cite{L4}, where a number of new 1-variable
equations in that range were given).

In our thesis \cite{G1} we gave the first (non-trivial) functional 
equations for $m=6$ (in fact, a whole family of equations in two variables) and for $m=7$, 
as well as many new examples for $\,m\le5\,$ (some of those results had already been included in \cite{Z1}
and in Chapter 16 of \cite{L2}). The main tool in our---computer-aided---investigation 
was Zagier's criterion for functional equations of the associated one-valued versions of the $m$-logarithm
(cf. \cite{Z1}, Proposition 1 and Proposition 3).

The relation of polylogarithms and algebraic $K-$theory, and in particular
Zagier's conjecture on special values of Dedekind zeta functions, have made it
clear that the question of finding functional equations for higher polylogarithms---as well as to 
relate them to each other---is a central one (cf.,~e.g.,~\cite{Z2}). 
For $m=3$, Goncharov \cite{Gon} found a basic functional equation as the key step in his proof of 
Zagier's conjecture for this case, while Wojtkowiak \cite{Wo} gave a whole family of relations.
The relationship between these equations has not been clarified so far.

The contents of the paper are as follows: in \S2, we recall a general functional equation for the 
dilogarithm, which is found to have some ``companion'' for the trilogarithm (Proposition \ref{wojtsymm})
which in turn specializes to Wojtkowiak's trilogarithm equations (Corollary \ref{wojt_eq}). We restate
(in \S3.1) Goncharov's equation in increasingly symmetric ways and relate it (in \S3.2 and \S3.3) to Wojtkowiak's
equa\-tions---the combining link being a 34-term equation, in fact a very special case of Wojt\-ko\-wiak's equation 
which turns out to be equivalent to Goncharov's one. In \S4, we state a family of 4-logarithm equations, a uniform 
and rather simple proof of which is given at the end of the paper (\S6). Finally, in \S5, we present our
``highscore result'' (obtained in 1992): a functional equation for the 7-logarithm in two variables (with 274 terms).

\section{A general dilogarithm equation}

\subsection{} Recall that the function  $\,\Li_m(z)\,$ is a many-valued function on 
$\,\BC-\{0,1\}\,$ but that one has the one-valued continuous function 
${\CL}_m:\C \to \R$ given by
$${\CL}_m(z)=\Re_m\bigg(\sum_{r=0}^m{ \frac{2^r B_r }{ r!}  \log^r|z| \,\Li_{m-r}(z) \bigg)},\qquad |z|\le 1,\ z\not\in\{0,1\},$$
where $\,\Re_m\,$ denotes the real part for $\,m\,$ odd and the imaginary part 
for $\,m\,$ even, and $\,B_r\,$ the $\,r$-th Bernoulli number 
($\,B_0=1,B_1=-1/2,B_2=1/6,\dots\,$). 
For $\,|z|> 1\,$, $\,{\CL}_m(z)\,$ is given by the functional equation
${\CL}_m(z)=(-1)^{m-1} {\CL}_m({1/z})\,, $ while for $z\in\{0,1,\infty\}$ one 
extends the function by continuity.
\smallskip
For $\,m=2\,$, this is the famous
Bloch--Wigner dilogarithm (cf., e.g., \cite{Bloch}). For $\,m>2\,$, one-valued versions of $\,\Li_m\,$
were introduced by Ramakrishnan,  Wojtkowiak and Zagier. The above function ${\CL}_m(z)$
was introduced by Zagier \cite{Z1} who denoted it $\,P_m(z)\,$. It has the advantage that
the corresponding functional equations of $\,\Li_m\,$ become ``clean'' for  $\CL_m$, 
e.g.,~in the 5-term relation stated above the elementary term on the right hand side disappears
when we replace $\,\Li_2\,$ by $\,{\CL}_2\,$.

Note that we can consider $\,{\CL}_m\,$ as a function on $\,\BC[\BC]\,$, the formal $\BC$-linear combinations
of elements in $\BC$, by extending it linearly.

\subsection{} There is a general functional equation for the dilogarithm found by Rogers ---for a
polynomial $\phi$ without constant term---and generalized by Wojtkowiak \cite{Wo} and Zagier
\cite{Z-D} to any rational function $\phi$ (cf.~also  Wechsung \cite{We1}, \cite{We2}). For $m\in \BN$ put
$$\widetilde{\L_m}(x,y,z,w)=\L_m\big(\crr(x,y,z,w)\big)=
\L_m\bigg(\frac{x-z}{ x-w}\,\frac{y-w}{ y-z}\bigg)\ ,$$
where $\crr$ denotes the cross ratio.
\begin{thm} Let $\,\phi:\BP^1_\BC\to\BP^1_\BC\,$ be a
rational function, $\,\alpha,B,C,D\in \BP^1(\BC)\,$. \\ Then,
denoting ${\rm deg}(\phi)$ the degree of $\phi$ and putting $\, A=\phi(\alpha)\,,$ we obtain
$$\sum_{\beta, \gamma, \delta}
{\widetilde{\L_2}(\alpha,
\beta,\gamma,\delta)}={\rm deg}(\phi)\cdot \widetilde{\L_2}(A,B,C,D)\,,$$ 
where $\beta$, $\gamma$ and $\delta$ run through the preimages of $B$, $C$ and $D$, respectively, with multiplicities.
\end{thm}
The proposition contains as a special case the 5-term relation (e.g.,~if we take $\phi(z)=z(1-z)$ 
and $B=1$, $C=0$, $D=\infty$) and in fact many other known equations.
A functional equation which is not covered by the proposition (and presumably the only one, essentially)
is the one relating $z$ and $\Bar{z}$, its complex conjugate, via $\L_2(\Bar{z})=-\L_2(z)$. We expect that the latter,
together with the equation in the above proposition, generates all functional equations for the dilogarithm.
In fact, a result by Wojtkowiak \cite{WoN} states that all functional equations for the dilogarithm with 
arguments $\C$-rational expressions in finitely many variables can be written as a sum of 5-term relations.

\section{Trilogarithm equations}
The following equation for the trilogarithm was found by symmetrizing 
a functional equation given by Wojtkowiak \cite{Wo}. 

 \begin{thm} \label{wojtsymm} Let $\,\phi:\BP^1_\BC\to\BP^1_\BC\,$ be a
rational function, $\,A_i,B_j,C_k,D_l\in\BP^1(\BC)\,$, $\,i,j,k,l\in\{1,2\}\,$.
Then, denoting ${\rm deg}(\phi)$ the degree of $\phi$, we have
$$\sum_{i,j,k,l}{(-1)^{i+j+k+l}\bigg(\sum_{\alpha_i,\beta_j,\gamma_k,\delta_l}
{\widetilde{\L_3}(\alpha,\beta,\gamma,\delta)}-
{\rm deg}(\phi)\cdot \widetilde{\L_3}(A_i,B_j,C_k,D_l)\bigg)}=0\,,$$
where $\alpha_i$, $\beta_j$, $\gamma_k$ and $\delta_l$ run through the preimages of $A_i$, $B_j$, $C_k$ and $D_l$, respectively, with multiplicities.
\end{thm}

\medskip
Wojtkowiak's original equation \cite{Wo}, pp.226-227, (a related equation had been found earlier by Wechsung, cf. \cite{We2}, \S.4) is obtained by specializing the equation from the proposition above
and can be stated in a simpler form as follows:

\begin{cor} \label{wojt_eq}
\ Let $\,A,B,C\in \BP^1(\BC)\,$, $\,\phi:\BP_\BC^1\to
\BP_\BC^1\,$ a rational function of degree $\,n\,$ and $\,x\,$ an independent
variable.\\
Let $\,\phi^{-1}(A)=\{\alpha_i\}_{i=1}^n\,$, 
$\,\phi^{-1}(B)=\{\beta_i\}_{i=1}^n\,$ and 
$\,\phi^{-1}(C)=\{\gamma_i\}_{i=1}^n\,$. Then
\begin{align*}
\widetilde{\L_3}\left(\phi(x),C,B,A\right)&-\sum_{i,j,k=1}^n {\widetilde{\L_3}\left(x,\gamma_i,\beta_j,\alpha_k\right)}\\
&-\sum_{i,j,k=1}^n {\widetilde{\L_3}\left(x,\alpha_i,\alpha_j,\gamma_k\right)}-
\sum_{i,j,k=1}^n {\widetilde{\L_3}\left(x,\beta_i,\beta_j,\gamma_k\right)}\\
&+\sum_{i,j,k=1}^n {\widetilde{\L_3}\left(x,\alpha_i,\alpha_j,\beta_k\right)}
+\sum_{i,j,k=1}^n {\widetilde{\L_3}\left(x,\beta_i,\beta_j,\alpha_k\right)}\\
=&\ {\rm const},
\end{align*}
i.e.,~the expression on the left hand side is independent of $\,x\,$.
\end{cor}
\newline
Here it is understood that $\,\alpha_i\,$, $\,\beta_j\,$ and $\,\gamma_k\,$ run through
all preimages (counted with multiplicity) of $\,A\,$, $\,B\,$ and 
$\,C\,$, respectively.\\

\subsection{Around Goncharov's equation}\label{around_gon}
\subsubsection{The original description.}
Goncharov \cite{Gon} found a beautiful interpretation of certain functional equations in terms
of configuration spaces, and as a crucial by-product he provided an equation for the 
trilogarithm in three variables $\,\al_{i}\,$ for which he gave a threefold 
symmetry. We reproduce it here using the shorthand 
$\be_i=1-\al_i+\al_i \al_{i-1}$ $(i=1,2,3)$, where indices are understood modulo 3: form
the formal linear combination $\gamma(\al_1,\al_2,\al_3)\in \Z[\BQ(\al_1,\al_2,\al_3)]$ given by
\begin{align}\label{gonrel}
\gamma(\al_1,\al_2,\al_3)=&{\sum_{i=1}^3} \Big( \Big[\frac{1}{ \al_i}\Big]+ [ \be_i] + 
\Big[ \frac{\al_i\al_{i-1}}{\be_i }\Big] + 
\Big[\frac{\be_i}{ \be_{i+1}\al_{i+2}} \Big] +\Big[-\frac{\be_i \al_{i+1}}{ 
\be_{i+1}}\Big]\Big)  \nonumber \\ 
&+\Big[-\frac{1}{ \al_1\al_2\al_3}\Big] -\sum_{i=1}^3\Big( \Big[\frac{\be_i}{ \al_{i-1}} \Big] + \Big[\frac{\be_i}{ \be_{i+1} \al_{i}\al_{i-1}}\Big] +[1]
\Big) \,.
\end{align}
\begin{prop} With the above notation, we have
 $$\L_3\big(\gamma(\al_1,\al_2,\al_3)\big)=0\,.$$
\end{prop}
Since there are 22 non-constant terms occurring in $\gamma(\al_1,\al_2,\al_3)$, we will refer to it as 
22-term (or Goncharov's) relation.
\subsubsection{A more symmetric description.}\label{moresym}
Since Goncharov's equation plays such a central role for the theory, it seems worthwhile to 
analyze its structure a bit further.\\
There is actually a much bigger symmetry group $G$ (of order 192) than the cyclic 
group on three letters acting on the set of arguments and dividing the 22
non-constant terms into 2 orbits, one of length 16 (corresponding to the 15 terms in the first sum 
in (\ref{gonrel}) plus the single term $\big[-\frac{1}{ \al_1\al_2\al_3}\big]$),
the other one of length 6
(corresponding to the second sum, with the exclusion of the constant terms $[1]$).

$G$ is generated by two involutions  $$\pi_1:(\al_1,\al_2,\al_3) 
\mapsto \Big(\al_1,\al_2,-\frac{\be_1}{ \al_1\be_3}\Big)\,,\qquad  \pi_2:(\al_1,\al_2,\al_3) \mapsto \Big(\frac{1}{ \al_1},\frac{1}{ \al_3},\frac{1}{
\al_2}\Big)\,,$$
together with the obvious symmetry of order 3 (shifting the indices mod 3).
 The set consisting of the arguments 
in the sixteen terms of the first orbit
is transformed either into itself (e.g.,~via $\pi_1$) or into the set of inverses 
(e.g.,~via $\pi_2$). We give a  presentation with more obvious symmetries:
let $\, t_1,\dots,t_4\,$ be four variables, subject to the constraint $\,\prod_i{}t_i=1\,$.
Then the following element in $\,\Q(t_1,t_2,t_3)\,$ is annihilated by $\L_3$ for each (meaningful) 
evaluation of $t_1,t_2,t_3$ in $\C$ 
\big(with $t_4=(t_1t_2t_3)^{-1}$\big):
$$ \sum_i{}[t_i]\,+\,\sum_{\scriptstyle{i,j\atop i\ne j}}{}\bigg[\frac{1-t_i}{ 1-t_j^{-1}}\bigg]\,-\,\frac{1}{ 4} \sum_{
\scriptstyle{i,j\atop i\ne j}}{}[t_i\,t_j]
\,-\,\frac{1}{ 8}\sum_{\{i,j,k,l\}=\{1,2,3,4\}}{}\bigg[\frac{(1-t_i)(1-t_j)}{ (1-t_k^{-1})(1-t_l
^{-1})}
\bigg]\,-3[1]\,. $$
The obvious $\,\Sig_4\,$-action on the set of arguments (where we identify 
$[z]$ and $[1/z]$),
together with the involution (cf.~\S1.1)
$$t_i\mapsto \frac{1-t_i}{ 1-t_{i+2}^{-1}}\, \qquad (i\ {\rm mod}\ 4\,)$$
generates the symmetry group $G$ (of order 192) mentioned above.

\subsubsection{A yet more symmetric description.}\label{yetmore}
Consider the finite group $G'$ (of order 96) of automorphisms of $\Q(y_1,y_2,y_3,z_1,z_2,z_3)$,
generated by $G'=\langle g,h\rangle$, where
\begin{align} \label{generators}
g:(y_1,y_2,y_3,z_1,z_2,z_3)&\mapsto(\frac{1}{ y_1},z_2,z_3,z_1,y_2,y_3),\ \qquad \\
\label{merkh}
h:(y_1,y_2,y_3,z_1,z_2,z_3)&\mapsto(y_2,y_3,y_1,z_2,z_3,z_1)\,.
\end{align}
The orbits of $\ y_1\ $ and of $\ \displaystyle{\vphantom{\Bigg|}\frac{y_1-z_3}{1-y_1z_2}\frac{y_2-z_1}{1-y_2z_3}\frac{y_3-z_2}{1-y_3 z_1}}\ $
under $\ G'\ $ have order 12 and 32, respectively.
As in the previous subsection, we introduce ``parametrization variables'' $t_1,\dots,t_4,$ subject to the constraint
$\prod_{j=1}^4 t_j=1$, and form 6 arguments
$$\{A_i\}_{i=1,2,3}=\{t_i t_4\}_i\,,\quad \{B_i\}_{i=1,2,3}=\Big\{\frac{1-t_j^{-1}}{1-t_i}
\frac{1-t_k^{-1}}{1-t_4}\,\Big|\, \{i,j,k\}=\{1,2,3\}\Big\}\,.$$
The involutory automorphism induced by 
$$(t_1,t_2,t_3,t_4)\mapsto \big(\iota(t_4,t_1),\iota(t_3,t_2),\iota(t_2,t_3),\iota(t_1,t_4)\big)$$
where $\displaystyle{\iota(x,y)=\frac{1-x}{1-y^{-1}}}$, acts on $(B_1,B_2,B_3,A_1,A_2,A_3)$ like the automorphism $g$ in 
(\ref{generators}) above 
while the automorphism induced by $(t_1,t_2,t_3)\mapsto (t_2,t_3,t_1)$ acts like $h$ in (\ref{merkh}).\\
If we put $y_i$ and $z_i$ equal to $A_i$ and $B_i$, respectively, in the above, then the $G'$-orbit of the expression $\ \displaystyle{\vphantom{\Bigg|}\frac{y_1-z_3}{1-y_1z_2}\frac{y_2-z_1}{1-y_2z_3}\frac{y_3-z_2}{1-y_3 z_1}}\ $ contains
only 16 different terms up to inversion of the arguments---obviously, the $G'$-orbit of $\ y_1\ $ contains only 6 terms 
up to inversion---and the resulting 16+6 arguments coincide precisely with the ones occurring in Goncharov's 
equation.
\subsubsection{Other descriptions.}
A different way to use the above parametrization with $A_i$ and $B_i$ in order to exhibit the symmetries among the 22 terms is obtained by putting $\alpha_i=-\sqrt{A_i}$, $\beta_i=\sqrt{B_i}$. Then the four elements 
$\alpha_1^{\eps_1} \alpha_2^{\eps_2} \alpha_3^{\eps_3}$ with $\eps_j\in\{\pm1\}$ $(j=1,2,3)$ and $\,\prod \eps_j =-1\,$, together with 
the twelve elements $\,\alpha_i^{\eps_i} \beta_{i+1}^{\eps_{i+1}} \beta_{i+2}^{\eps_{i+2}}\,$ \\ ($i \mod 3$) 
with $\eps_j\in\{\pm1\}$ and $\prod \eps_j =1$, 
form the 16 arguments in the second orbit above in \S\ref{yetmore}, while the $\alpha_i^2$ and $\beta_i^2$ give the 6 arguments in the
first orbit.

Yet another description---in terms of equations instead of parametrizations---is given by the following 9 (redundant)
equations where we put $\displaystyle{q(x,y)=\frac{x-y}{1-xy}}$:
\begin{align*}
&B_i^{-1}=\frac{q(B_{i-1},B_{i+1})}{q(A_{i+1},A_{i-1})}\,,\\
&A_i^{-1}=q(A_{i+1},\frac1 {B_{i-1}})\ q(A_{i-1},\frac 1{B_{i+1}})\,,\\
&\frac{A_{i+1}(1-A_i)^2}{A_i(1-A_{i+1})^2} =\frac{B_{i+1}(1-B_i)^2}{B_i(1-B_{i+1})^2}\,,
\end{align*}
where $i=1,2,3$ and the indices are taken modulo 3. (Obviously, we can break the symmetry and
reduce this to a  system of 7 equations in only 4 variables by eliminating, say, 
$A_1$ and $B_1$. The system of equations is no complete intersection, though.)

\subsection{Relating Goncharov's and Wojtkowiak's equations}\label{gontowojt}
In order to compare the equations resulting from Goncharov's and Wojtkowiak's approach, we propose
to study an ``intermediate'' relation, with 34 terms, which allows an interpretation for both situations.

1. In Wojtkowiak's equation, the ``intermediate'' relation occurs when we consider
 the rational function
$$\phi(x)=\frac{x-a}{x-c^{-1}}\cdot \frac{x-b}{x-abc}=\,\crr(x,xbc,b,abc)\,,$$
(for $x=1$ this is a $B_i$-type term above) and $(A,B,C)=(\infty,0,1)$, and subtracting two different specializations for $x$
of the ensuing terms from Corollary \ref{wojt_eq}---each specialization providing 17 terms---we are left with 34 terms which form a functional equation in 
three variables. The generic 17 terms are given, in factored form, as the arguments of the following formal
linear combination in $\Z[F]$, $F=\BQ(a,b,c,t)$:
\newcommand \wrap[1] {\Big[#1\Big]}
\begin{align*}
&f(a,b,c,t) =\wrap {\frac{(1-ct)a}{a-t}} + \wrap {\frac{(1-ct)b}{b-t}}+\wrap {\frac{1-ct}{c(a-t)}} + \wrap {\frac{1-ct}{c(b-t)}} \nonumber \\
&\phantom{f(a,b,c,t) =} +\wrap{\frac{abc-t}{(a-t)bc}} +\wrap{\frac{abc-t}{(b-t)ac}} +\wrap{
\frac{abc-t}{a-t}} +\wrap{\frac{abc-t}{b-t}} \\+&\wrap{\frac{(a-t)b(ac-1)}{(b-t)a(bc-1)}} +\wrap{\frac{(t-a)(1-bc)}{(t-b)(1-ac)}}
+\wrap{\frac{(tc-1)b(ac-1)}{(abc-t)(bc-1)}} +\wrap{\frac{(tc-1)a(bc-1)}{(abc-t)(ac-1)}}  \nonumber\\
&\ -\wrap{\frac{(1-ct)a b c}{a b c -t}} -\wrap{\frac{(1-ct)}{(abc-t)c}} -\wrap{\frac{t-a}{t-b}} -\wrap{\frac{b(a-t)}{a(b-t)}} 
-\wrap{\frac{(b-t)(a-t)c}{(abc-t)(1-ct)}}\,. \nonumber
\end{align*}
\ms
\begin{df} The {\em 34-term relation} is given by the difference $$f(a,b,c,t)-f(a,b,c,u) \in \Z[F]\,,\qquad 
F=\BQ(a,b,c,t,u).$$
\end{df}
\ms
\begin{rem} From Corollary \ref{wojt_eq} it results that $\CL_3$ vanishes on 34-term relations.
\end{rem}

2. Let us recall that Goncharov \cite{Gon} has defined a triple ratio for 6 points in the projective plane 
over a field $F$, with values in $\Z[F]$, in such a way
that it respects a 7-term relation: for any 7 distinct points $P_1,\dots,P_7$ in $\PP^2$ one has
$$\sum_{i=1}^7 (-1)^i\cdot \text{triple ratio}(P_1,\dots,\widehat{P_i},\dots,P_7)=0\,.$$
In the setting of configurations, the 34-term relation mentioned above encodes the well-definedness of the triple ratio
associated to a configuration of six points $P_i$ ($i=1,\dots,6$) in general position in the plane as follows.

Goncharov reduces such a configuration with the help of the 7-term relation above to more degenerate configurations
by introducing the intersection point $Q$ of the line through the points $P_1,P_2$ with the line
passing through $P_3$ and $P_4$. He uses the fact that he has already defined the triple ratio for 
the more degenerate configurations which result by leaving out any other of the 6 original points $P_i$.
But there are different possibilities for choosing $Q$, depending on the ordering of the $P_i$---e.g.,~one could 
switch the roles of $P_2$ and $P_3$---and one verifies (this is not explicitly done in \cite{Gon}) that the difference of two such
possibilities, if non-zero, essentially results in two sets of $17$ terms which correspond precisely to the 34 terms in question.
(Here ``essentially'' alludes to the fact that we argue modulo two simple functional equations for the
trilogarithm, the inversion relation  $\L_3(x)=\L_3(1/x)$
and the 3-term relation 
$\L_3(x)+\L_3\big(1/(1-x)\big)+\L_3(1-1/x)=\L_3(1)$.)

\subsection{Relating the 34-term and 22-term equation} In this subsection, we will indicate a proof of the following fact:
\begin{prop} The 22-term relation and the 34-term 
relation are equivalent in the sense that each one, together with its specializations, implies the other.
\end{prop}

If we take the difference of the specializations of $f(a,b,c,t)$ to $t=1$ and $t=0$, respectively, then a number
of terms degenerate and we are left with the sum of a Kummer-Spence relation and a 22-term relation. (Again, we
work up to inversion and 3-term relation, the latter being a specialization of both 22-term and 34-term equation.) Similarly, the difference of the specializations to $t=a$ and $t=0$ gives
a version of the Kummer-Spence relation itself. Thus the 34-term relation implies the 22-term relation. 
Conversely, the substitution 
$$(t_1,t_2,t_3)\mapsto \big(\crr(t,0,c^{-1},a),\crr(t,0,b,c^{-1}),\crr(t,0,a,abc)\big)$$ maps
the terms of the 22-term relation above in the form given in \S\ref{moresym} to the 17 terms of $f(a,b,c,t)$ 
together with five terms (plus a constant term) which are independent of $t$. 
This shows immediately that the 22-term relation implies the 34-term relation.

Combining the proposition with the previous subsection, we conclude:

\begin{prop} Goncharov's 22-term relation is subsumed in Wojtkowiak's family of equations in Corollary \ref{wojt_eq}.
\end{prop}

\subsection{A 21-term equation}
Yet another related equation of possible interest is a sum of four 22-term relations which combines to 
a functional equation in 3 variables with less, viz.~only 21, different non-constant arguments (up to inversion), but with 
coefficients in $\{\pm1, \pm2\}$:
with $\gamma$ as in (\ref{gonrel}), consider
$$\Gamma(x,y,z)=\gamma\Big(\frac1 {1-x},\frac{1-x}{1-x y},1-z\Big)+\gamma\Big(1-\frac 1 x,\frac{1-xy}{y(1-x)},\frac 1 {1-z^{-1}}\Big)\,$$
and note that its symmetrization in the first 2 arguments, which of course also constitutes a functional equation for the
trilogarithm, has the following form---we introduce a ``symmetrizing variable''
$z_2=(x_1x_2z_1)^{-1}$ and put $\ j(t,u)=\frac{1-u^{-1}}{1-t}\ $:
\begin{align*}
\Gamma(x_1,x_2,z_1)+&\Gamma(x_2,x_1,z_1)=\ -2[x_1\,x_2] -2[1]\\
+2\sum_{i=1}^2\Big(&[x_i]+[j(x_i,x_{3-i})]+[z_i]+[j(z_i,z_{3-i})]\\
-&[x_i\,j(z_1,z_2)]-[j(x_i,x_{3-i})\,j(z_1,z_2)]-[x_iz_1]-[j(x_i,x_{3-i})\,z_1]\Big)\\
+\sum_{i=1}^2 \Big(&[x_i\,z_1\,j(x_i,x_{3-i})\,j(z_1,z_2)] + [ j(x_i,x_{3-i})\,j(z_1,z_2)\,/\,(x_i\,z_1)]\Big)\,.
\end{align*}

\section{4-logarithm equations}
It is in general a tedious job to verify even a single functional equation for the polylogarithms. Therefore
it seems worthwhile finding a rather short way to verify a whole family of them.
We propose to do this for a family of equations for the 4-logarithm, using essentially only one
(and actually trivially checked) polynomial equation. One can proceed similarly (albeit in a somewhat more
complicated fashion) with other families (cf. \cite{Ga}, Thm.~4.4 and Thm.~4.9), even up to $n=6$.
Let us emphasize that not only had there been no example for $n\geq6$ previously known 
(apart from the trivial distribution relations) but Wechsung \cite{We1} even {\sl proved} that
the type of equation (called ``Kummer-type'') which gave the only examples for $n=4, 5$ known at the time is not ``good enough''
for  $n\geq6$. The following example, which is proved in \S6, is not of Kummer-type.

\begin{thm}\label{fourlog} Let $\,F\,$ be a field of characteristic 0, $\phi(x)=
x^{n-1}(x-1)$ for some
$n\in \BN$.
For some $t,u\in F$, let 
$\,\{x_i\}_{i=1}^{n}=\phi^{-1}(t)\,$ resp. $\,\{y_i\}_{i=1}^{n}=\phi^{-1}(u)\,$ be the sets of
preimages in some finite extension field $F'$ of $F$. 
Then the following element in $\BZ[F']$ is annihilated by ${\CL}_4$:
\begin{multline}
n(n-2)\bigg[{\prod_i{}x_i\over \prod_j{}y_j}\bigg] -(n-1)^2
\sum_{i,j=1}^{n}\bigg[{1-x_i^{-1}\over 1-y_j^{-1}}\bigg]
+n^2 \sum_{i,j=1}^{n}\bigg[{1-x_i\over 1-y_j}\bigg]\\ \label{thmeq}
-n^2(n-1)^2 \sum_{i,j=1}^{n}\bigg[{x_i\over y_j}\bigg]
+ n(n-1)^2\sum_{i=1}^{n}{}\Big(\big[1-{1\over x_i}\big]-\big[1-{1\over 
y_i}\big]\Big)
\,. \end{multline}
\end{thm}

\section{A 7-logarithm equation}

In this section we give a 2-variable functional equation for the 7--logarithm. 
This functional equation consists of 274 terms, but by using an appropriately 
symmetrized notation we can write it in a more digestible form.

\subsection{}\label{sevenlog} Let $\,F\,$ be the rational function field in two variables $\,\BQ(t,u)\,$.
Put

$$ f_1(z)=\frac{-z}{ 1-z+z^2},\quad f_2(z)=\frac{z-1}{ 1-z+z^2}\,,\quad
f_3(z)=\frac{z(1-z)}{ 1-z+z^2}\,$$
$${\rm and}\quad \,f(z)=-f_1(z)f_2(z) f_3(z)=\frac{z^2(1-z)^2}{ (1-z+z^2)^3}\,.$$ 
Notice that the group $G\cong \Sig_3$ generated by $z\to 1/z$ and $z\to \,1-z$ 
permutes the $f_i$ and leaves $f$ invariant.
For $\,a,b,c,d\in\BZ\,$ define elements $\,\li{a,b;c,d}\}_0(t,u)\ $ and 
$\,\li{a,b;c,d}\}\,$ in $\,\BC[\BC]\,$ by 
$$\,\li{a,b;c,d}\}_0(t,u)=\sum_{i,j=1}^3{}
{\bigg[\frac{f(t)^a f_i(t)^{b-a}
}{ f(u)^c f_j(u)^{d-c}}\bigg]}\quad {\rm and}$$
\begin{align*}
\,\li{a,b;c,d}\}&=\li{a,b;c,d}\}_0(t,u)+\li{c,d;a,b}\}_0(t,u)\\
&=\li{a,b;c,d}\}_0(t,u)+\li{a,b;c,d}\}_0(u,t)\,,\end{align*}
where the latter equation holds only up to inversions of the arguments (which, of course,
does not affect the functional equation to be presented below).
\smallskip\noindent
Finally, set
\def \bl{{\phantom{\frac{1}{ 2}}}}

\halign{\hfill{}$\strut#$&\hfill$#$&\hfill$#$&\hfill$#$&\hfill$#$&\hfill$#$&
\hfill$#$&\hfill$#$&\hfill$#$&\hfill$#$&\hfill$#$&\hfill$#$&\hfill$#$&
\hfill$#$&\hfill$#$&\hfill$#$&\hfill$#$&\hfill$#$&\hfill$#$&\hfill$#$&\hfill$#$&\hfill$#$&$#$\hfill\cr
\,\xi_{7}^{(3)}\,=\,-\,\frac{1}{ 18}\,\cdot&\frac{609}{ 4}& \li&-1,&-1;&-1,&-1\re& -\,\frac{1}{ 3}\,\cdot&35& 
\li& -1,&-1;&-2,&1\re &
+\,\frac{1}{ 3}\,\cdot&\frac{105}{ 8}& \li& -1,&-1;&3,&-5 \re& \cr  
-\,\frac{1}{ 3}\,\cdot&21& \li& -1,&-1;&-1,&4\re&
-\,\frac{1}{ 3}\,\cdot&15& \li& -1,&-1;&-2,&5\re& +\,\frac{1}{ 3}\,\cdot&15& \li& -1,&-1;&3,&-4\re,&\cr
\cr
\xi_{7}^{(2)}\,=\,+\,\frac{1}{ 2}\,\cdot&700& \li& 1,&0;&1,&0\re& +\,\frac{1}{ 2}\,\cdot&\frac{175}{ 4}& \li& 1,&-3;&
1,&-3\re&
+\,\frac{1}{ 2}\,\cdot&28& \li& -2,&3;&-2,&3\re& \cr
-\bl&35& \li& 1,&-3;&-2,&3\re& -\bl&140& \li& -2,&3;&1,&0\re&
+\bl &175& \li& 1,&0;&1,&-3\re ,&\cr
\cr
\xi_{7}^{(1)}\,=\,+\,\frac{1}{ 2}\,\cdot&700& \li& 1,&-2;&-1,&2\re& +\bl &3150& \li& 0,&1;&1,&-1\re& 
+\,\frac{1}{ 2}\,\cdot&\,1575& \li& -1,&1;&1,&-1\re&\cr
-\bl&2100& \li& 1,&-2;&0,&-1\re& +\,\frac{1}{ 2}\,\cdot&\,6300& \li& 0,&1;&0,&-1\re& 
-\bl&\,1050& \li& -1,&2;&-1,&1\re&\cr
-\,\frac{1}{ 2}\,\cdot&700& \li& -1,&2;&-1,&2\re &-\,\frac{1}{ 2}\,\cdot&
\,1575 &\li& -1,&1;&-1,&1\re& -\,\frac{1}{ 2}\,\cdot&\,6300& \li& 0,&1;&0,&1\re&\cr
+\bl &\,1050 &\li& -1,&2;&1,&-1\re& +\bl &2100& \li& 0,&-1;&-1,&2\re& -\bl&3150 &\li& 1,&-1;&0,&-1\re&.\cr }
\smallskip\noindent
 
\begin{thm}\label{crown} The element $\,\xi_{7}=\xi_7^{(1)}+\xi_7^{(2)}+\xi_7^{(3)}\,\in \BQ[F]$ is a functional equation for
$\,{\CL}_7\,$, i.e.,~
$${\CL}_7\big(\xi_{7}(t,u)\big)=0$$ 
for any $\,t,u\in \BC\,$ where $\,f_j(t),f_j(u)\notin\{0,\infty\}\,$, $\,j=1,2,3\,$.
\end{thm}

\begin{proof}
Functional equations for polylogarithms can be characterized by an algebraic
criterion (cf.~\cite{Z1}, Proposition 1), so the proof can be reduced to checking
this criterion for $\,\xi_7\,$ which was done with a computer program.
The statement applies, of course, to any functional equation; the difficulty
lies in {\em finding} functional equations, not in checking their validity.
\end{proof}
 

\subsection{Remark} \ \ In $\,\xi_7\,$ the denominator of the first factor of a 
coefficient (which we suppress if it is 1) 
gives the multiplicity of an argument in the sum $\,\{a,b;c,d\}\,$
(where we identify $\,[z]\,$ and 
$\,[1/z]\,$).

\newline
We associate a weight to each $\,\{a,b;c,d\}\,$ by 
$\,\,{\rm wt}(\{a,b;c,d\})={\rm wt}(a,b)\,$
in the following way:  setting
$\,{\rm wt}(a,b)=\frac{1}{ 2}(|a|+|a+b|+|2a+b|)\,$ we find that in all terms  
$\,\{a,b;c,d\}\,$ of $\,\xi_{7}\,$ we have $\,{\rm wt}(a,b)={\rm wt}(c,d)\,$.
\newline
This provides each $\xi_7^{(i)}$ ($i=1,2,3$) with weight $i$.
There is a weight-preserving operation on these blocks which also includes the coefficients.
We explain this operation in more detail in the following subsection.\\

\subsection{More symmetric
 form of $\,\xi_{7}\,$.\ }

\medskip
Let $\,(\BZ^3)_0=\{(a,b,c)\in\BZ^3\mid a+b+c=0\}\,$ and 
$\,(\BZ^3)_1=\{(a,a,b)\in\BZ^3\}\,$. We consider the map 
\begin{align}\Theta:\BZ^3&\to \BZ^3\,, \\
(a,b,c)&\mapsto (a,-b-c,b-a) \end{align}
which sends $\,(\BZ^3)_0\,$ to $\,(\BZ^3)_1\,$,
and a second one (the $\,f_i$ are taken as in (\ref{sevenlog}))
\begin{align}\Phi:\BZ^3&\to \{\phi:\BP^1_\BC\to\BP^1_\BC\mid \phi\ {\rm rational}\}\,, \\
\alpha=(\alpha_1,\alpha_2,\alpha_3)&\mapsto \Big(\phi_\alpha:z\mapsto \big(-f_1(z)\big)
^{\alpha_1} f_2(z)^{\alpha_2} f_3(z)^{\alpha_3}\Big)\,. \end{align}
\newline
The symmetric group $\,\Sig_3\,$ operates on $\,(\BZ^3)_0\,$ and $\,\BZ^3\,$ via permutation. 
\newline
For each $\,k\in\BZ\,$ we have the element $\,(k,-1,1-k)\in(\BZ^3)_0\,$.
\newline We set
$$A_k=\Theta\big(\Sig_3\cdot (k,-1,1-k)\big)\,,\qquad k=1,2,3.$$
\newline
Then we have
\begin{align*} A_1&=\{\pm(1,1,-2),\pm(-1,-1,1),\pm(0,0,1)\}\,, \\
A_2&=\{(2,2,-3),(-1,-1,3),(-1,-1,0)\}\,,\\
A_3&=\{(-1,-1,-1),(-1,-1,4),(-2,-2,5),(-2,-2,1),(3,3,-4),(3,3,-5)\}\,.
\end{align*}
We will denote by $\delta$ the $\Sig_3$-invariant element $(-1,-1,-1)$ of $A_3$ which will play a special role.
Finally we define for $\,\alpha=(\alpha_1,\alpha_2,\alpha_3)\not=0\,$
a ``weight''
$$\omega(\alpha)=\frac{1}{ \alpha_1-\alpha_3}\,.$$
Then Theorem \ref{crown} can be restated more concisely by saying that $$\ 60\,\xi_{7}=\sum_{i=1}^3 
\xi^{(i)}\ $$
with
\begin{align*}
\xi^{(1)}=&\ -\frac{29}{ 20}\Big[\frac{\phi_\delta(t)}{ \phi_\delta(u)}\Big]
+\sum_{\alpha\in A_3\setminus\{\delta\}}{{\omega(\alpha)}}\sum_{\sigma\in \Sig_3/\Sig_2}{}
\Big(\Big[
\frac{\phi_{\sigma\alpha}(t)}{ \phi_\delta(u)}\Big]+\Big[\frac{\phi_\delta(t)}{ \phi_{\sigma\alpha}(u)}\Big]\Big)
\,,\cr
\xi^{(2)}= &\,\phantom{-}\,\, \frac{20}{ 3}\sum_{\alpha,\beta\in A_2}{}{ \omega(\alpha)
\omega(\beta)}\sum_{\sigma,\tau\in \Sig_3/\Sig_2}{}
\Big[\frac{\phi_{\sigma\alpha}(t)}{ \phi_{\tau\beta}(u)}\Big]\,,\cr
\xi^{(3)}=&\,-\, 30\sum_{\alpha,\beta\in A_1}{}{ \omega(\alpha)
\omega(\beta)}\sum_{\sigma,\tau\in \Sig_3/\Sig_2}{}
\Big[\frac{\phi_{\sigma\alpha}(t)}{ \phi_{\tau\beta}(u)}
\Big]\,.\end{align*}

\noindent
Altogether we get $\,1+30+81+162=274\,$ different arguments (up to inverses) since in the
last sum (over $\,A_1\,$)  
the arguments associated to $\,(\alpha,\beta)\,$ and $\,(-\alpha,-\beta)\,$ are inverse
to each other.

\section{A proof of Theorem \ref{fourlog}}

{\bf Proof.}  \ Using Zagier's criterion (cf.~\cite{Z1}, Proposition 1) alluded to earlier, we need only
verify that the above combination lies in the kernel of the map
\begin{align*} \beta_4:\ \BZ[F']&\longrightarrow \Sym^2 {F'}^\times\otimes\,\bigwedge
{}^2\, 
{F'}^\times\\
[x]&\longmapsto x^{\odot 2}\otimes \big(x\wedge (1-x)\big)\,.\end{align*}
Here the $\bigwedge{}^2$ denotes the second exterior product (i.e., two-fold tensors subject to the relations
$x\wedge x=0$). Note that $xy\wedge z = x\wedge z + y\wedge z$.

Our strategy is as follows: a few preliminary considerations (steps 0--2) allow to rewrite the $\beta_4$-images in a more 
convenient way, so that the theorem is essentially reduced to Claim \ref{claim} (in step 3).

{\bf Part I: Reformulations.}  0. Note that $t-\phi(z)$ is a polynomial in (the variable) $z$ with roots
equal to $\{x_i\}_i$, therefore 
$$t-\phi(z)=\phi(x_i)-\phi(z) = \lambda \prod_i (x_i-z)$$
for some constant $\lambda$. Thus, for any fixed $l$ and $m$, 
$$ t-u=\phi(x_l)-\phi(y_m) =  \lambda \prod_i (x_i-y_m)=  \mu \prod_j (x_l-y_j) $$
for some constants $\lambda$ and $\mu$ (which actually turn out to be equal to $\pm 1$ and therefore can be
neglected in the following).

1. As a first preliminary step, we express $\alpha$ and $1-\alpha$ in terms of the factors $x_i$
and $y_j$ for each of the arguments $\alpha$ in (\ref{thmeq}). In order to save indices, we put
$x=x_l$ and $y=y_m$ for some fixed $l,m \in \{1,\dots,n\}$. We obtain
$$\frac{1-x}{1-y} = \frac{y^{n-1}X}{x^{n-1}Y}\qquad \text{and}\qquad \frac{1-x^{-1}}{1-y^{-1}} = \frac{y^{n}X}{x^{n}Y}\,,$$
where we have set $X=\prod_i x_i$, $Y=\prod_j y_j$. A further decomposition we need is 
$$1-\frac{1-x}{1-y} =\frac{y^{n-1}}{Y}(x-y) \qquad \text{and}\qquad 1-\frac{1-x^{-1}}{1-y^{-1}} = \frac{y^{n-1}}{xY}(x-y) \,.$$
With these preparations, we can write the corresponding images under $\beta_4$ as
\begin{equation}\label{erste}
\beta_4\bigg(\Big[\frac{1-x}{1-y}\Big]\bigg)= \left(\frac{y^{n-1}X}{x^{n-1}Y}\right)^{\otimes 2}\otimes
\frac{y^{n-1}X}{x^{n-1}Y}\wedge  \frac{y^{n-1}}{Y}(x-y)
\end{equation}
and
\begin{equation}\label{zweite}
\beta_4\bigg(\Big[\frac{1-x^{-1}}{1-y^{-1}}\Big]\bigg)= \left(\frac{y^{n}X}{x^{n}Y} \right)^{\otimes 2}\otimes
\frac{y^{n}X}{x^{n}Y} \wedge  \frac{y^{n-1}}{xY}(x-y)\,.
\end{equation}
These two expressions decompose ``naturally''  into two parts, one of which contains the factor $(x-y)$ in the
last tensor factor. E.g., the part in $\Wedge^2{F'}^\times $ in the first expression (\ref{erste}) factors as
\begin{equation}\label{dritte}
 \frac{y^{n-1}X}{x^{n-1}Y}\wedge  \frac{y^{n-1}}{Y}(x-y) =  \frac{X}{x^{n-1}}\wedge  \frac{y^{n-1}}{Y}
\ +\ \frac{y^{n-1}X}{x^{n-1}Y}\wedge  (x-y)\,
\end{equation}
where the first summand on the right has been reduced using the defining property of $\,\wedge$. 

2. As a second preliminary step, we ``switch'' to an additive notation---formally, we put $\xi_i= ``\log x_i$'', 
$\eta_j= ``\log y_j$'', and introduce $\zeta_{lm}= ``\log(x_l-y_m)$'' as well as the shorthands $\xi=\sum_i\xi_i$
and $\eta=\sum_j\eta_j$ (note also that the $\wedge$ now satisfies $(a+b)\wedge c = a\wedge c+ b\wedge c$); 
then the term (\ref{erste}) resp. (\ref{zweite}) is written, using the decomposition (\ref{dritte}), as
\begin{multline}\label{vierte}
\big((n-1)(\eta_m-\xi_l)+\xi -\eta\big)^3 \wedge \zeta_{lm}\\
 + \big((n-1)(\eta_m-\xi_l)+\xi -\eta\big)^2
\cdot \Big(\xi - (n-1)\xi_l\Big)\wedge\Big(-\eta +(n-1)\eta_m\Big)
\end{multline}
resp.
\begin{multline}\label{fuenfte}
\big(n(\eta_m-\xi_l)+\xi - \eta\big)^3 \wedge \zeta_{lm}\\
 + \big(n(\eta_m-\xi_l)+\xi - \eta\big)^2
\cdot \big(\xi - (n-1)\xi_l +\eta_m\big)\wedge\big(-\eta+(n-1)\eta_m- \xi_l\big)\,.
\end{multline}
Introducing the shorthands $\ S=\xi-\eta\ $ and $\ s_{lm} = \xi_l-\eta_m\ $, we rewrite the latter two 
equations as
\begin{equation}\label{sechste}
\big(S- (n-1)s_{lm}\big)^3 \wedge \zeta_{lm} + \big(S- (n-1)s_{lm}\big)^2
\cdot \big(\xi - (n-1)\xi_l\big)\wedge\big(-\eta+(n-1)\eta_m\big)
\end{equation}
resp.
\begin{equation}\label{siebte}
\big(S- ns_{lm}\big)^3 \wedge \zeta_{lm} + \big(S- ns_{lm}\big)^2
\cdot \big( \xi - (n-1)\xi_l+\eta_m\big)\wedge\big(-\eta+(n-1)\eta_m-\xi_l\big)\,.
\end{equation}

{\bf Part II: Calculations.} 3. So far, we have only reformulated the objects under consideration. Now we can proceed to the two 
actual calculations involved, the first of which follows easily from  the readily verified identity
\begin{equation} \label{achte}
\big(\xi-(n-1)\xi_l\big) \wedge \big(\eta -(n-1)\eta_m\big) = \sum_{i,j} (2-n)^{\delta_{il}+\delta_{jm}} 
\xi_i\wedge \eta_j\,.
\end{equation}
Here $\delta_{ij} (=1\text{ if }i=j,\text{ and } =0\text{ otherwise})$ denotes the usual ``Kronecker-$\delta$''.

Using relation (\ref{achte}) together with the presentations (\ref{sechste}) and (\ref{siebte}), 
the $\beta_4$-image of the first two terms in (\ref{thmeq})
becomes
\begin{equation}\label{beta_vier}
\beta_4\bigg(n^2\Big[\frac{1-x_l}{1-y_m}\Big] -(n-1)^2 \Big[\frac{1-x_l^{-1}}{1-y_m^{-1}}\Big]\bigg)
=T_1^{lm}+T_2^{lm} +T_3^{lm} +T_4^{lm}  
\end{equation}
where 
\begin{eqnarray*}
T_1^{lm}&=&\big((2n-1)S^3 - 3n(n-1)S^2 s_{lm} +n^2(n-1)^2s_{lm}^3\big)\wedge \zeta_{lm}\,,  \\
T_2^{lm}&=& -\Big(n^2\big(S-(n-1)s_{lm}\big)^2 - (n-1)^2 (S-ns_{lm})^2\Big)\bigg(\sum_{i,j}(2-n)^{\delta_{il}+\delta_{jm}}\ \xi_i\wedge\eta_j
\bigg)\,,\\
T_3^{lm}&=& (n-1)^2 (S-ns_{lm})^2\cdot \big(\eta_m\wedge \xi_l \big)\,,  \\
T_4^{lm}&=&  (n-1)^2 (S-ns_{lm})^2\cdot \big(\xi\wedge \xi_l\ +\ \eta_m\wedge  \eta\big)\,.
\end{eqnarray*}

We will sum the latter expression over all $l$ and $m$ (unless explicitly indicated otherwise, 
all sums in the following run from 1 to $n$).
\begin{claim} \label{claim} With notations as above, we have
\begin{align*}
\sum_{l,m}\beta_4\bigg(n^2\Big[\frac{1-x_l}{1-y_m}\Big]& -(n-1)^2 \Big[\frac{1-x_l^{-1}}{1-y_m^{-1}}\Big]\bigg)\\
 =&\ 
n(n-2)\Bigg(-S^3\wedge \Big(\frac 1 n \sum_{l,m}\zeta_{lm}\Big) + S^2\sum_{l,m} \xi_l\wedge \eta_m\Bigg)\\
&+n^2(n-1)^2\Bigg(\sum_{l,m}s_{lm}^3\wedge \zeta_{lm}\ -\ \sum_{l,m}s_{lm}^2(\xi_l\wedge \eta_m)\Bigg) \\ 
& + \sum_i f_1(\xi_i) + \sum_j f_2(\eta_j)
\end{align*}
for some functions $f_k$ ($k=1,2$) depending on only one of the variables $t$ and $u$.
\end{claim}
\begin{proof} We will need the following immediately verified identities:
\begin{equation} \label{ident1}
Z:=\frac 1 n \sum_{i,j}\zeta_{ij} = \sum_i \zeta_{im} =\sum_j\zeta_{lj} \,,\quad\forall l\,,\  \forall m\,,
\end{equation}
\begin{equation} \label{ident2}
S= \sum_{i}\xi_{i} -\sum_j\eta_j = \frac 1 n \sum_{i,j} s_{ij} \,,
\end{equation}
\begin{equation} \label{ident0}
\sum_{m}(\xi-n\xi_{m})=\sum_l(\eta-n\eta_l) =0\,,
\end{equation}
and the ``distribution properties'' (for fixed $l$, $m$):
\begin{equation} \label{ident3}
\sum_{i,j} (2-n)^{\delta_{il}+\delta_{jm}} = \sum_i (2-n)^{\delta_{il}} =1\,,
\end{equation}
as well as some simple consequence of the above
\begin{equation} \label{ident4}
\sum_{l,m} s_{l m}(2-n)^{\delta_{il}+\delta_{jm}} = S - (n-1)s_{ij}\,.
\end{equation}
[Proof: write $s_{lm}=\xi_l-\eta_m$ and use successively the equations in (\ref{ident3}).]

We will first treat the sum $\sum_{l,m} T_4^{lm}$ which is easily reduced to the form 
$ \sum_i f_1(\xi_i) + \sum_j f_2(\eta_j)$ since each of its terms which is dependent on both $\xi$'s and $\eta$'s
is cancelled in the sum by virtue of (\ref{ident0}).

Thus it will suffice to identify the sum $\sum_{i=1}^3 \sum_{l,m}T_i^{lm}$ with the first two lines on the
right hand side of the claim. The first two terms of
\begin{multline} \label{T1}
\sum_{l,m}T_1^{lm} = (2n-1)S^3\wedge \sum_{l,m} \zeta_{lm} - 3n (n-1)S^2 \sum_{l,m} s_{lm}\wedge \zeta_{lm} \\
+n^2(n-1)^2
\sum_{l,m} s_{lm}^3\wedge \zeta_{lm}
\end{multline}
combine to 
\begin{equation*}
-n(n-2) S^3\wedge Z
\end{equation*}
\big(note that $\sum_{l,m} s_{lm}\wedge \zeta_{lm} = S\wedge Z$ by using (\ref{ident1}) and (\ref{ident2})\big),
while the third term on the right of (\ref{T1}) is already of the desired form.

The term 
\begin{equation*}
\sum_{l,m}T_2^{lm} = -\sum_{l,m}\bigg((2n-1)S^2-2n(n-1)Ss_{lm}\bigg)\sum_{i,j}(2-n)^{\delta_{il}+\delta_{jm}}\ \xi_i\wedge \eta_j
\end{equation*}
decomposes into two simpler sums, the first of which is given as
\begin{equation}\label{quadratic}
-(2n-1)S^2\sum_{i,j}\sum_{l,m}(2-n)^{\delta_{il}+\delta_{jm}}\ \xi_i\wedge \eta_j\,=-(2n-1)S^2\sum_{i,j} \xi_i\wedge \eta_j\,,
\end{equation}
where we have used (\ref{ident3}), while the second one equals
\begin{multline} \label{linear}
2n(n-1)S \sum_{l,m}s_{lm}\sum_{i,j}(2-n)^{\delta_{il}+\delta_{jm}}\ \xi_i\wedge \eta_j=\\
2n(n-1)S^2 \sum_{i,j} \xi_i\wedge \eta_j\ - \ 2n(n-1)^2S\sum_{i,j} s_{ij}\big(\xi_i\wedge \eta_j\big)
\end{multline}
by  virtue of (\ref{ident4}). Finally, expanding $\sum_{l,m}T_3^{lm}$ as a ``polynomial'' in $S$
gives us three further sums, one being quadratic, one linear and one constant in $S$; 
the quadratic one combines with (\ref{quadratic}) and with the first
term on the right in (\ref{linear}) to $\ n(n-2)S^2\sum_{i,j}\xi_i\wedge \eta_j\ $, the linear one cancels with the second term in
(\ref{linear}), and the constant one is given by $\ - n^2(n-1)^2 \sum_{l,m} s_{lm}^2 (\xi_l\wedge \eta_m)\ $. 
This proves Claim \ref{claim}.

4. We finish by observing that the first two terms on the right hand side of Claim \ref{claim} correspond precisely
to the images under $\beta_4$ of $-n(n-2)\big[\prod x_i/\prod y_j\big]$ and $n^2(n-1)^2\sum_{i,j}\big[x_i/y_j\big]$, respectively,
while $\sum_{l,m}T_4^{lm}$, which corresponds to the last two terms of the claim, can be recognized readily---after more careful investigation of the terms ``purely in $\xi_i$ or $\eta_j$''---
as the image under $\beta_4$ of $ -n(n-1)^2\big(\sum_i\big[1-x_i^{-1}\big] - \sum_j\big[1-y_j^{-1}\big]\big)$ .

\end{proof}

\noindent{\Small {\bf Acknowledgements}: We would like to heartily thank Don Zagier for not only suggesting
the original problem, but even more for his steady support, encourage- and enlightenment, and constant flow of ideas 
surrounding this work. It was made possible through the generosity and hospitality of the 
Max-Planck-Institut f\"ur Mathematik in Bonn with its unsurpassed working conditions.}

\bibliographystyle{plain}

\end{document}